\documentclass[11pt,reqno]{amsart}
\usepackage{amsmath,amssymb,amsfonts,amscd,latexsym,amsthm,mathrsfs}
\usepackage[usenames]{color}
\usepackage[unicode]{hyperref}
\renewcommand{\d}{{\mathrm d}}

\newcommand{\ZZ}{{\mathbb Z}}

\newtheorem{theorem}{Theorem}

\allowdisplaybreaks[2]

\begin{document}

\title{A case study for $\zeta(4)$}

\author{Carsten Schneider}
\address{Johannes Kepler University Linz, Research Institute for Symbolic Computation, Altenberger Str. 69, A-4040 Linz, Austria}
\email{carsten.schneider@risc.jku.at}

\author{Wadim Zudilin}
\address{Department of Mathematics, IMAPP, Radboud University, PO Box 9010, 6500~GL Nij\-megen, Netherlands}
\email{w.zudilin@math.ru.nl}

\date{17 April 2020. \emph{Revised}: 22 September 2020}

\thanks{Research partially supported by the Austrian Science Fund (FWF) grants SFB F5006-N15 and F5009-N15 in the framework of the Special Research Program ``Algorithmic and Enumerative Combinatorics''.}

\begin{abstract}
Using symbolic summation tools in the setting of difference rings,
we prove a two-parametric identity that relates rational approximations to~$\zeta(4)$.
\end{abstract}

\maketitle

\hbox to\hsize{\hfill\vbox{\hsize=64mm\vskip1mm
		\noindent
		\hbox to 20mm{Kingdom:\hfill}Mathematical constants\\
		\hbox to 20mm{Class:\hfill}Periods\\
		\hbox to 20mm{Family:\hfill}Multiple zeta values\\
		\hbox to 20mm{Genus:\hfill}Single zeta values\\
		\hbox to 20mm{Species:\hfill}Even zeta values}}

\section{Introduction}
\label{sec1}

The quantity
$$
\zeta(4)=\sum_{k=1}^\infty\frac1{k^4}=\frac{\pi^4}{90}
$$
is a somewhat typical representative of even zeta values\,---\,the values of Riemann's zeta function at positive even integers.
It is shadowed by the far more famous $\zeta(2)=\pi^2/6$, which was a main subject of Euler's resolution of Basel's problem,
and $\zeta(3)$\,---\,an \emph{objet de l'\'etude} of Ap\'ery's iconic proof of the irrationality of the latter (and also of $\zeta(2)$) \cite{Ap79,vdP79}.
Though known to be irrational (and transcendental!), $\zeta(4)$ serves as a natural guinea pig for extending Ap\'ery's machinery to other zeta values.
Ap\'ery-type approximations to the number were discovered and rediscovered on several occasions \cite{Co81,So02,Zu03}, however they were not good enough to draw conclusions about its irrationality.
An unexpected difficulty to control the `true' arithmetic of those rational approximations to $\zeta(4)$ generated further research \cite{KR07,Zu09}, which eventually led to producing sufficient approximations and establishing a new world record for the irrationality measure of $\pi^4$ \cite{MZ20}.

In this note we turn our attention to a rational side of the coin and prove the following two-parametric identity.

\begin{theorem}
\label{th:main}
For integers $n\ge m\ge0$, define two rational functions
\begin{align}
R(t)=R_{n,m}(t)
&=(-1)^m\Bigl(t+\frac n2\Bigr)\frac{(t-n)_m}{m!}\,\frac{(t-2n+m)_{2n-m}}{(2n-m)!}\,\nonumber
\\ &\qquad\times
\frac{(t+n+1)_n}{(t)_{n+1}}\,\frac{(t+n+1)_{2n-m}}{(t)_{2n-m+1}}\,\biggl(\frac{n!}{(t)_{n+1}}\biggr)^2\nonumber
\\ \intertext{and}
\tilde R(t)=\tilde R_{n,m}(t)
&=\frac{n!\,(t-n)_{2n-m}}{(t)_{n+1}(t)_{2n-m+1}}
\sum_{j=0}^n{\binom nj}^2\binom{2n-m+j}n\frac{(t-j)_n}{n!}\,.\label{Equ:RTilde}
\end{align}
Then
\begin{equation}
-\frac13\sum_{\nu=n-m+1}^\infty\frac{\d R(t)}{\d t}\bigg|_{t=\nu}
=\frac16\sum_{\nu=1}^\infty\frac{\d^2\tilde R(t)}{\d t^2}\bigg|_{t=\nu}.
\label{hyp1}
\end{equation}
\end{theorem}

The $m=n$ instance of \eqref{hyp1} was stated as Problem~1 in \cite{Zu09}.

The fact that both sides of \eqref{hyp1} are linear forms in $1$ and $\zeta(4)$ with rational coefficients is verifiable by standard techniques \cite{KR07,Zu03,Zu09} which employ the partial-fraction decomposition of the rational functions.
A remarkable outcome of this identity is the \emph{coincidence} of two different-looking rational approximations to the zeta value.
Such coincidences are often a source of deep algorithmic and analytical developments\,---\,check \cite{EZZ20} for another exploration of this theme (see also \cite{BCS20}).

The main difficulty in establishing equality \eqref{hyp1} (in contrast to tackling, for example, Ap\'ery's sums  in~\cite{AperyCA} for $\zeta(3)$) is that its both sides are not hypergeometric functions but rather \emph{derivatives} of hypergeometric functions. Another issue is that the summation range on the left-hand side is somewhat unnatural.

\section{Symbolic summation}
\label{sec2}

Denote by $Z_l(n,m)$ and $Z_r(n,m)$ the left- and right-hand sides of~\eqref{hyp1}, respectively. 
In order to prove the identity~\eqref{hyp1} we proceed as follows.

\smallskip
\noindent
\textbf{(A)} We compute the linear recurrence
\begin{equation}\label{Equ:ZRec}
a_0(n,m) Z(n,m)+a_1(n,m) Z(n,m+1)+a_2(n,m) Z(n,m+2) = 0
\end{equation}
with
\begin{equation}\label{Equ:ZRecCoeff}
\begin{split}
a_0(n,m)&=(2n-m)^5,
\\
a_1(n,m)&=-(4n-2m-1) (6n^4-24n^3m+22n^2m^2-8nm^3+m^4-24n^3
\\ &\qquad
+30n^2m-14nm^2+2m^3+8n^2-10nm+2m^2-4n+m),
\\
a_2(n,m)&=-(2n-m-1)^3 (4n-m) (m+2),
\end{split}
\end{equation}
which holds simultaneously for $Z(n,m)=Z_l(n,m)$ and $Z(n,m)=Z_r(n,m)$ for all $n,m\in\ZZ_{\geq0}$ with $n-2\ge m\ge0$. In addition, we observe that $a_2(n,m)\neq0$ for all $n,m\in\ZZ_{\geq0}$ with $0\leq m<n$.

\smallskip
\noindent
\textbf{(B)} We show that the following initial values hold:
\begin{align}
Z_l(n,0)&=Z_r(n,0) \quad\text{for all}\; n\geq0,
\label{Equ:Initial0}\\
Z_l(n,1)&=Z_r(n,1) \quad\text{for all}\; n\geq1.
\label{Equ:Initial1}
\end{align}
Combined with \textbf{(A)} this proves that 
$Z_l(n,m)=Z_r(n,m)$
holds true for all $n\ge m\ge0$.

\smallskip
In order to carry out the steps \textbf{(A)} and \textbf{(B)}, advanced symbolic summation techniques in the setting of difference rings are utilized. Among them the following three summation paradigms play a decisive role, that are available within the summation package~\texttt{Sigma}~\cite{Schneider:07a}.

\noindent
\textbf{(i) Creative telescoping.}
Given a sum $F(m)=\sum_{\nu=a}^bf(m,\nu)$ and $\delta\in\ZZ_{\geq0}$, one searches for polynomials $c_0(m),\dots,c_{\delta}(m)$, free of $\nu$, and $g(m,\nu)$ such that
\begin{equation}\label{Equ:SummandRecurrence}
g(m,\nu+1)-g(m,\nu)=c_0(n)f(m,\nu)+c_1(m)f(m+1,\nu)+\dots+c_{\delta}(m)f(m+\delta,\nu)
\end{equation}
holds for all $a\leq\nu\leq b$.
Thus summing~\eqref{Equ:SummandRecurrence} over $\nu$ one obtains the recurrence
\begin{equation}\label{Equ:SumRecurrence}
g(m,b+1)-g(m,a)=c_0(m)F(m)+c_1(m)F(m+1)+\dots+c_{\delta}(m)F(m+\delta).
\end{equation}
By specializing $a,b$ further\,---\,e.g., to $a=0$ and $b=m$, or sending $b$ to $\infty$ if the limit exists\,---\,one obtains recurrence relations for more specific sums. The computed creative telescoping solution $(c_0(m),\dots,c_{\delta}(m),g(m,\nu))$ is also called a proof certificate for the recurrence~\eqref{Equ:SumRecurrence} found: usually it allows one to verify that $F(m)$ is a solution of~\eqref{Equ:SumRecurrence} by simple polynomial arithmetic, without analyzing the usually complicated computation steps of the underlying summation algorithm. The algorithmic version of creative telescoping has been introduced in~\cite{Zeilberger:91,AequalB} for hypergeometric sums. In order to prove~\eqref{hyp1}, we will employ a generalized machinery for creative telescoping~\cite{Schneider:15} where the summand can be composed not only in terms of hypergeometric products, but of indefinite nested sums defined over hypergeometric products. We emphasize that all recurrences produced below (using the \texttt{Sigma}-command \texttt{GenerateRecurrence}) are accompanied by such proof certificates which guarantee the correctness of all the calculations. Since the output is rather large and can be easily reproduced with \texttt{Sigma}, any explicit printout of the proof certificates is skipped. 

\noindent
\textbf{(ii) Recurrence solving.}
Given a linear recurrence of the form~\eqref{Equ:SumRecurrence}, one can search for solutions that are expressible within certain classes function spaces. Using the \texttt{Sigma}-command \texttt{SolveRecurrence} one can search for hypergeometric solutions~\cite{Petkov:92,AequalB} and, more generally, for all solutions that are expressible in terms of indefinite nested sums defined over hypergeometric products. Such solutions are also called d'Alembertian solutions~\cite{Abramov:94,Schneider:01} a subclass of Liouvillian solutions~\cite{Singer:99}. 

\noindent
\textbf{(iii) Simplification of expressions.}
Within \texttt{Sigma} the expressions in terms of indefinite nested sums defined over hypergeometric products are represented in the setting of difference rings and fields~\cite{Karr:81,Schneider:08c,DR1}. 
Utilizing this difference ring machinery~\cite{Schneider:10c,DR3} (compare also~\cite{Singer:08}) one can apply, e.g.,
the \texttt{Sigma}-command \texttt{SigmaReduce} to an expression in terms of indefinite nested sums. Then the output is a simplified expression where the arising sums and products (except products such as $(-1)^m$) are independent among each other as functions of their external parameter. In particular, the input expression evaluates to zero (from a certain point on) if and only if \texttt{Sigma} reduces the expression to the zero-expression.

\smallskip
These summation paradigms can be used to transform a definite (multi-)sum to an expression in terms of indefinite nested sums by deriving a linear recurrence, solving the recurrence found in terms of indefinite nested sums, and, in case that sufficiently many solutions are found, combining them to an expression that evaluates to the same sequence as the input sum. Recently this machinery has been used for large scale problems coming from particle physics (see, e.g.,~\cite{CALadder:16} and references therein). In this regard, also the package \texttt{EvaluateMultiSum}~\cite{Schneider:13a}, which automatizes this summation mechanism, has been utilized non-trivially in the sections below.

In the following sections we present the main steps of our proof for Theorem~\ref{th:main} that is based on the above summation algorithms. All the necessary calculation steps are collected in a Mathematica notebook that can be accessed via%
\footnote{In case that the reader does not have access to Mathematica, we supplement the pdf file \href{https://www.risc.jku.at/people/cschneid/data/SchneiderZudilinMMA.pdf}{\texttt{SchneiderZudilinMMA.pdf}} (same www-path!) that contains all the calculations in printed form.}
\begin{center}
	\url{https://www.risc.jku.at/people/cschneid/data/SchneiderZudilinMMA.nb}\,.
\end{center}

\section{A linear recurrence in $m$ for the left-hand side}
\label{Subsec:Z_l}

In order to activate the summation package \texttt{Sigma}, the sums arising in~\eqref{hyp1} have to be tailored to an appropriate input format. As it turns out below, one can carry out the differentiation by introducing additionally the harmonic numbers
$$
S_a(n)=\sum_{k=1}^n\frac1{k^a}
$$
of order $a\in\ZZ_{\ge0}$. 
Though we see no natural way to obtain such a representation for the full summation range $\nu$ with $n-m+1\leq\nu$, splitting it into the ranges over $\nu$ with $n-m+1\leq\nu\leq 2n-m-1$ and $2n-m\leq\nu$ makes the job well.
More precisely, we split the left-hand side of~\eqref{hyp1} into the two subsums
\begin{align*}
W_1(n,m)&=\sum_{\nu=2n-m+1}^{\infty}\frac{\d R_{n,m}(t)}{\d t}\bigg|_{t=\nu}=\sum_{\nu=1}^{\infty}\frac{\d R_{n,m}(t+2n-m)}{\d t}\bigg|_{t=\nu}
\\ \intertext{and}
W_2(n,m)&=\sum_{\nu=n-m+1}^{2n-m}\frac{\d R_{n,m}(t)}{\d t}\bigg|_{t=\nu}=\sum_{\nu=1}^{n}\frac{\d R_{n,m}(t+n-m)}{\d t}\bigg|_{t=\nu},
\end{align*}
so that
\begin{equation}\label{ZW1Ws}
Z_l(n,m)=-\frac1{3}\big(W_1(n,m)+W_2(n,m)\big).
\end{equation}
Observe that 
\begin{multline*}
R_{n,m}(t+2n-m)=
(-1)^m\Bigl(t+2n-m+\frac n2\Bigr)\frac{(t+n-m)_m}{m!}\,\frac{(t)_{2n-m}}{(2n-m)!}
\\
\times
\frac{(t+3n-m+1)_n}{(t+2n-m)_{n+1}}\,\frac{(t+3n-m+1)_{2n-m}}{(t)_{2n-m+1}}\,\biggl(\frac{n!}{(t+2n-m)_{n+1}}\biggr)^2
\end{multline*}
and 
\begin{multline*}
R_{n,m}(t+n-m)=(-1)^m\Bigl(t+n-m+\frac n2\Bigr)\frac{(t-m)_m}{m!}\,\frac{(t-n)_{2n-m}}{(2n-m)!}
\\ \qquad\times
\frac{(t+2n-m+1)_n}{(t+n-m)_{n+1}}\,\frac{(t+2n-m+1)_{2n-m}}{(t+n-m)_{2n-m+1}}\,\biggl(\frac{n!}{(t+n-m)_{n+1}}\biggr)^2.
\end{multline*}
By definition all Pochhammer symbols in the former expression are of the form $(t+x)_k$ for some $x\in\ZZ_{>0}$ and $k\geq0$. Thus, we can apply the formula
\begin{equation}\label{Equ:DPochhammer}
\frac{\d}{\d t}(x+t)_k\big|_{t=\nu}=(x+\nu)_k\big(S_1(\nu+x+k-1)-S_1(\nu+x-1)\big)
\end{equation}
for $\nu\in\ZZ$ with $x+\nu\in\ZZ_{>0}$ which follows from the product-rule of differentiation.
Employing this formula we get for all $\nu=1,2,\dots$ the following representation:  
\begin{align*}
&
F_1(n,m,\nu)
=\frac{\d}{\d t}R_{n,m}(t+2n-m)\bigg|_{t=\nu}
\\ &\quad
=\frac{(-1)^m n!^2 (1+\nu )_{-1	-m+2 n} (-m+n	+\nu)_m (1-m+3 n+\nu)_n (1-m+3 n+\nu)_{-m+2 n}}
{2\,m! (-m+2 n)!(-m+2 n+\nu )_{1+n}^3 (-m+2 n+\nu)_{1	-m+2 n}}
\\ &\quad\;\times
\bigg(-6 \nu
+\nu  (-2 m+5 n+2 \nu)
\big(
-S_1({\nu })
-S_1({-m+n+\nu })
+5 S_1({-m+2 n+\nu })
\\ &\quad\;\quad
-5 S_1({-m+3 n+\nu })
-S_1({-2 m+4 n+\nu })
+S_1({n+\nu })
+S_1({-m+4 n+\nu })
\\ &\quad\;\quad
+S_1({-2 m+5 n+\nu })
\big)
+\frac{5 n (m-2 n)}{m-2 n-\nu}
+\frac{n (-2 m+3 n)}{n+\nu}
+\frac{3 n (m-n)}{-m+n+\nu}
\bigg).
\end{align*}
Further, we prepare the summand of $W_2(n,m)$. Notice that the rule~\eqref{Equ:DPochhammer} cannot be applied to the arising factor $(t-n)_{2n-m}$. However we can easily overcome this issue by using the following elementary identity:
For $\nu\in\ZZ_{>0}$ with $1\leq \nu\leq n$ and any differentiable function $f(t)$, we have
\begin{equation}\label{Equ:SpecialD1}
\frac{\d}{\d t}\big((t-n)_{2n-m}f(t)\big)\bigg|_{t=\nu}
=(-1)^{n-\nu}f(\nu)(\nu+n-m-1)!(n-\nu)!\,.
\end{equation}
Therefore, for all $\nu\in\ZZ_{>0}$ with $1\leq \nu\leq n$ we get
\begin{align*}
F_2(n,m,\nu)&=\frac{\d R_{n,m}(t+n-m)}{\d t}\\
&=(-1)^m\Bigl(\nu+n-m+\frac n2\Bigr)\frac{(\nu-m)_m}{m!}\,\frac{(-1)^{n-\nu}(\nu+n-m-1)!(n-\nu)!}{(2n-m)!}\\ 
&\quad\times
\frac{(\nu+2n-m+1)_n}{(\nu+n-m)_{n+1}}\,\frac{(\nu+2n-m+1)_{2n-m}}{(\nu+n-m)_{2n-m+1}}\,\biggl(\frac{n!}{(\nu+n-m)_{n+1}}\biggr)^2.
\end{align*}
Because of the factor $(\nu-m)_m$, we have $F_2(\nu)=0$ for all $\nu\in\ZZ_{>0}$ with $1\leq \nu\leq m$.
Consequently, $W_1(n,m)$ and $W_2(n,m)$ can be written as
\begin{equation*}
W_1(n,m)=\sum_{\nu=1}^{\infty}F_1(\nu)
\quad\text{and}\quad
W_2(n,m)=\sum_{\nu=m+1}^{n}F_2(\nu)=\sum_{\nu=1}^{n-m}F_2(\nu+m),
\end{equation*}
where the summands $F_1(\nu)$ and $F_2(\nu)$ are given in terms of hypergeometric products and linear combinations of harmonic numbers.
Since these sums fit the input class of \texttt{Sigma}, we can apply the command \texttt{GenerateRecurrence} to both sums and compute for $0\leq m\leq n$ the recurrences
\begin{equation*}
a_0(n,m) W_s(n,m) + a_1(n,m) W_s(n,m+1)+a_2(n,m) W_s(n,m+2) = r_s(n,m)
\quad\text{for}\; s=1,2,
\end{equation*}
where the coefficients are given in~\eqref{Equ:ZRecCoeff} and where $r_1(n,m)=-r_2(n,m)$ is too large to be reproduced here (verification of the latter equality required an extra simplification step with \texttt{Sigma}).
To compute the recurrence for the \emph{hypergeometric} sum $W_2(n,m)$ one can alternatively use the Mathematica package~\texttt{fastZeil} \cite{PauleSchorn:95} based on~\cite{Zeilberger:91}.

Thus, $Z_l(n,m)$ given in~\eqref{ZW1Ws} is a solution of the recurrence~\eqref{Equ:ZRec}.
For this part we needed 15~minutes to compute both recurrences and to combine them to~\eqref{Equ:ZRec}. 

\section{A linear recurrence in $m$ for the right-hand side}
\label{Subsec:Z_r}
In order to calculate a linear recurrence for $Z_r(n,m)$ we follow the same strategy as for $Z_l(n,m)$ in Section~\ref{Subsec:Z_l} by utilizing more advanced summation tools of \texttt{Sigma}. 
Collecting all products in~\eqref{Equ:RTilde} to
$$
G_{n,m,j}(t)=\frac{n!\,(t-n)_{2n-m}}{(t)_{n+1}(t)_{2n-m+1}}{\binom nj}^2\binom{2n-m+j}n\frac{(t-j)_n}{n!},
$$
the right-hand side of~\eqref{hyp1} can be rewritten as
\begin{align*}
Z_r(n,m)
&:=\frac16\sum_{\nu=1}^{\infty}\sum_{j=0}^{n}\frac{\d^2}{\d t^2}G_{n,m,j}(t)\bigg|_{t=\nu}.
\\
\intertext{Similarly to the previous section, we split the sum further into subsums (see~\eqref{Equ:ZrSumExpr} for the final split) such that the differential operator acting on the summands can be replaced by modified summands in terms of harmonic numbers. On the first step, we write}
Z_r(n,m)
&=\frac16\big(C_1(n,m)+C_2(n,m)\big)
\end{align*}
with
\begin{align*}
C_1(n,m)=\sum_{\nu=1}^{\infty}\sum_{j=0}^{n}\frac{\d^2}{\d t^2}G_{n,m,j}(t+n)\bigg|_{t=\nu}
\quad\text{and}\quad
C_2(n,m)=\sum_{\nu=1}^{n}\sum_{j=0}^{n}\frac{\d^2}{\d t^2}G_{n,m,j}(t)\bigg|_{t=\nu}
\end{align*}
and apply, as before, formula~\eqref{Equ:DPochhammer} and its relatives to get a monster summand of $C_1(n,m)$ (that fills two pages) in terms of the harmonic numbers of order $1$ and~$2$. For illustration we print out only a few lines:
\begin{align*}
G_1(n,m,j,\nu)
=&\frac{\d^2}{\d t^2}G_{n,m,j}(t+n)\bigg|_{t=\nu}
\\ 
=&\quad
\frac{\binom{n}{j}^2 \binom{j
		-m
		+2 n
	}{n} (\nu )_{-m
		+2 n
	} (-j
	+n
	+\nu
	)_n}{(n
	+\nu
	)_{1+n} (n
	+\nu
	)_{1
		-m
		+2 n
}}\bigg(\cdots\\
&+S_1({-j+n+\nu })^2
+S_1({-j+2 n+\nu })^2
+S_1({-m+2 n+\nu })^2\\
&+S_1({n+\nu })
\frac{4(-j^2 m n
	+2 j^2 n^2
	+\dots
	+m \nu ^3
	-7 n \nu ^3
	-2 \nu ^4)}{\nu  (n
	+\nu
	) (-j
	+n
	+\nu
	) (-j
	+2 n
	+\nu
	) (-m
	+2 n
	+\nu
	)}\\
&+\cdots\bigg).
\end{align*}
In order to tackle the summand of $C_2(n,m)$, we have to differentiate $G_{n,m,j}(t)$ twice.
With $p(t)=(t-n)_{2n-m}$ and 
\begin{equation}\label{Equ:qDef}
q(t)=\frac{G_{n,m}(t)}{p(t)}=\frac{n!}{(t)_{n+1}(t)_{2n-m+1}}{\binom nj}^2\binom{2n-m+j}n\frac{(t-j)_n}{n!}
\end{equation}
we conclude that for all $1\leq\nu\leq n$ we have
\begin{align*}
\tilde{G}(\nu)
&=\frac{\d^2}{\d t^2}G_{n,m,j}(t)\bigg|_{t=\nu}
=q(t)\,\frac{\d^2p(t)}{\d t^2}+2\frac{\d p(t)}{\d t}\,\frac{\d q(t)}{\d t}+p(t)\,\frac{\d^2q(t)}{\d t^2}\bigg|_{t=\nu}
\\
&=q(t)\frac{\d^2p(t)}{\d t^2}+2\frac{\d p(t)}{\d t}\frac{\d q(t)}{\d t}\bigg|_{t=\nu};
\end{align*}
the last equality follows since $p(t)|_{t=\nu}=0$ for all $1\leq\nu\leq n$.
Similarly to~\eqref{Equ:SpecialD1}, we can use in addition the following calculation:
For $\nu\in\ZZ_{>0}$ and $1\leq\nu\leq n$, we have
$$
\frac{\d}{\d t}(t-n)_{2n-m}\bigg|_{t=\nu}
=h(t)\bigg|_{t=\nu}
\quad\text{and}\quad
\frac12\,\frac{\d^2}{\d t^2}(t-n)_{2n-m}\bigg|_{t=\nu}
=\frac{\d}{\d t}h(t)\bigg|_{t=\nu}
$$
with 
$$h(t)=\frac{(-1)^{n-\nu}\Gamma(t+n-m)(\nu-t+1)_{n-\nu}}{\Gamma(t-\nu+1)}.$$
In particular, if $\nu>j$, we can apply the rule~\eqref{Equ:DPochhammer} to all Pochhammer symbols in~\eqref{Equ:qDef}: \begin{align*}
&
G_2(n,m,j,\nu)=\tilde{G}(\nu)
\\ &\quad
= 2q(t)\frac{\d}{\d t}h(t)+2h(t)\frac{\d}{\d t} q(t)\Big|_{t=\nu}
\\ &\quad
= \frac{2(-1)^{n+\nu } \binom{n}{j}^2 \binom{j-m+2 n}{n} (1)_{n-\nu} (2)_{-1-m+n+\nu} (1-j+\nu)_{-1+n}}
{\nu ^3 (-m+n+\nu)^2(1+\nu )_n (1+\nu )_{-m+2 n}}
\\ &\quad\;\times
\bigg(
\nu  (-j+\nu ) (-m+n+\nu ) \Big(\frac{1}{j-n-\nu}-S_1({-j+\nu })+S_1({-j+n+\nu })\Big)
\\ &\quad\quad
+\nu  (-j+\nu ) (-m+n+\nu ) \big(-S_1({-m+2 n+\nu })+S_1({\nu })\big)
\\ &\quad\quad
+\nu  (-j+\nu ) (-m+n+\nu ) \big(S_1({\nu })-S_1({n+\nu })\big)
\\ &\quad\quad
-\nu  (-j+\nu )+\nu  (-m+n+\nu )
\\ &\quad\quad
+2 (j-\nu ) (-m+n+\nu )
+\nu  (-j+\nu ) (-m+n+\nu )
\\ &\quad\quad
+\nu  (-j+\nu ) (-m+n+\nu ) \big(-1+S_1({-m+n+\nu })\big)
\\ &\quad\quad
-\nu  (-j+\nu ) (-m+n+\nu ) S_1({n-\nu })
\bigg).
\end{align*}
For $1\leq\nu\leq j$, we use $q(\nu)=0$ and apply the rule
$$
\frac{\d}{\d t}\bigl((t-j)_{n}f(t)\bigr)\bigg|_{t=\nu}
=f(\nu)(\nu-j)_{j-\nu} (n+\nu-j-1)!
$$
(compare with \eqref{Equ:SpecialD1})
valid for any differentiable function $f(t)$, in place of \eqref{Equ:DPochhammer}, to~\eqref{Equ:qDef}.
It follows that
\begin{align*}
G_3(n,m,j,\nu)=&\tilde{G}(\nu)={2\binom nj}^2\binom{2n-m+j}n\\
&\times\frac{(-1)^{n+\nu }(
	n
	+\nu-m-1
	)! (n
		-\nu
	)! (n
		+\nu-j-1
	)! (\nu-j)_{j
		-\nu
}}{(\nu )_{1+n} (\nu )_{
		2 n-m+1
}}.
\end{align*}
Therefore,
\begin{align*}
C_2(n,m)
&=\sum_{\nu=1}^{n}\sum_{j=0}^{n}\frac{\d^2}{\d t^2}G_{n,m}(t)\bigg|_{t=\nu}
\\
&=\sum_{j=0}^{n-1}\sum_{\nu=j+1}^{n}G_2(n,m,j,\nu)
+\sum_{j=1}^{n}\sum_{\nu=1}^{j}G_3(n,m,j,\nu),
\end{align*}
hence
\begin{align}
Z_r(n,m)
&=\frac16\big(C_1(n,m)+C_2(n,m)\big)
\nonumber\\
&=\frac16\bigg(\sum_{j=0}^n\sum_{\nu=1}^{\infty}G_1(n,m,j,\nu)
+\sum_{j=0}^{n-1}\sum_{\nu=j+1}^{n}G_2(n,m,j,\nu)
\nonumber\\ &\quad
+\sum_{j=1}^{n}\sum_{\nu=1}^{j}G_3(n,m,j,\nu)\bigg).
\label{Equ:ZrSumExpr}
\end{align}
Denote by $A_1(n,m)$, $A_2(n,m)$ and $A_3(n,m)$ the three resulting sums in~\eqref{Equ:ZrSumExpr}
and use \texttt{Sigma} to compute three linear recurrences of $A_{s}(n,m)$ with $s=1,2,3$.
A routine calculation demonstrates that each of the recurrences found can be brought to the form
\begin{equation}\label{Equ:ARecs}
a_0(n,m) A_s(n,m) + a_1(n,m) A_s(n,m+1)+ a_2(n,m) A_s(n,m+2) = u_s(n,m),
\end{equation}
where the coefficients are given in~\eqref{Equ:ZRecCoeff} and where only the right-hand sides $u_s(n,m)$ for $s=1,2,3$ differ.
As an illustration, we provide with details about how we treat
$$
A_1(n,m)=C_1(n,m)=\sum_{j=0}^n\sum_{\nu=1}^{\infty}G_1(n,m,j,\nu).
$$
In the first step, \texttt{Sigma} is used to compute a linear recurrence of the inner sum 
\begin{equation}\label{Equ:cDef}
c(n,m,j)=\sum_{\nu=1}^{\infty}G_1(n,m,j,\nu)
\end{equation}
in~$j$,
\begin{align}
&
(j-n)^2 (j-n+1)^2 (j-m+2n+1) (j-m+2n+2)c(n,m,j)
\nonumber\\ &\;
-(j-n+1)^2 (j-m+2n+2)
\big(2 j^3-2 j^2 m+2 j m n-3 j n^2+m n^2-2 n^3
\nonumber\\ &\; \qquad\qquad\qquad
+8 j^2-5 j m-2 j n+4 m n-7 n^2
+11 j-3 m-4 n
+5
\big)c(n,m,j+1)
\nonumber\\ &\;
+(j+2)^3(j-2n+1) (j-m+n+2)^2c(n,m,j+2)=r(n,m,j),
\label{Equ:RecPureInJ}
\end{align}
and one additional recurrence with one shift in $m$ and one shift in~$j$,
\begin{align}
&
(j-n)^2 (j-m+2n+1)
\big(
j^3+j m^2-j^2 m-m^3-2 j m n+4 m^2 n-4 m n^2
\nonumber\\ &\; \qquad\qquad\qquad
+2 j^2-j m-2 j n+2 m n-4 n^2
+j-2 n
\big)c(n,m,j)
\nonumber\\ &\;
-(j+1)^3(j-2 n) (j-m+n+1)^2c(n,m,j+1)
\nonumber\\ &\;
-(j-n)^2(j-m+2n) (j-m+2n+1) (m+1)(m-2n) c(n,m+1,j)
=s(n,m,j);
\label{Equ:RecmPureInJ}
\end{align}
here the right-hand sides $r(n,m,j)$ and $s(n,m,j)$ are large expressions in terms of hypergeometric products and the harmonic numbers 
$S_ 1({n})$, $S_ 1({2 n})$, $S_ 1(n-j)$, $S_ 1({2 n-j})$, $S_ 1({2 n-m})$, $S_1({3 n-m})$.
Finally, we use new algorithms that are described in~\cite{BRS:18} and that are built on ideas from~\cite{Schneider:05d,APS:05}.
Activating these new features of \texttt{Sigma} we can compute the linear recurrence~\eqref{Equ:ARecs} with $s=1$ where the right-hand side $u_0(n,m)$ is an expression in terms of the harmonic numbers $S_1({n}),S_1({2 n}),S_1({2 n-m}),S_1({3 n-m})$, the infinite sums 
\begin{equation}\label{Equ:InfiniteSums}
c(n,m,0), \; c(n,m,1), \; c(n,m,n+1)
\end{equation}
and the definite sums
\begin{equation}\label{Equ:FiniteSums}
\begin{aligned}
&\sum_{i=0}^n \binom{n}{i}^2 \binom{2n-m+i}{n} \frac{(n-i+1)_{n}}{(2n-i)^k}&&\quad\text{for}\; k=0,1,2,
\\
&\sum_{i=0}^n \binom{n}{i}^2 \binom{2n-m+i}{n} \frac{(n-i+1)_{n}}{(2n-i)^k}S_1(n-i)&&\quad\text{for}\; k=0,1,
\\
&\sum_{i=0}^n \binom{n}{i}^2 \binom{2n-m+i}{n} \frac{(n-i+1)_{n}}{(2n-i)^k}S_1(2n-i)&&\quad\text{for}\; k=0,1.
\end{aligned}
\end{equation}
Note that all these definite sums in~\eqref{Equ:FiniteSums} are \emph{not} expressible in terms of hypergeometric products and indefinite nested sums defined over such products. 
For example, the linear recurrence for the last sum in~\eqref{Equ:FiniteSums} with $k=0$ computed with \texttt{Sigma} has order $5$ and has not even a hypergeometric product solution.
We further remark that the above approach is connected to the classical holonomic summation approach~\cite{Zeilberger:90a} and their improvements given in~\cite{Chyzak:00,Koutschan:13}.
In all these traditional versions one needs systems composed by homogeneous recurrences. However, the transformations of~\eqref{Equ:RecPureInJ} and~\eqref{Equ:RecmPureInJ} to such a form would lead to gigantic recurrence systems and the computation of the desired linear recurrence~\eqref{Equ:ZRec} would be out of scope.

Using this refined holonomic summation approach with \texttt{Sigma}, we needed in total 10 minutes to derive the recurrence for $A_1(n,m)$ which holds for all $0\leq m\leq n$. Similarly, one can compute for the other two double sums $A_2(n,m)$ and $A_3(n,m)$ the recurrence~\eqref{Equ:ARecs} in 15 and 2 minutes, respectively, which hold for all $0\leq m\leq n-2$.
Here the right-hand sides $u_2(n,m)$, $u_3(n,m)$ consist of similar definite sums as given in~\eqref{Equ:FiniteSums}.
Adding up \eqref{Equ:ARecs} corresponding to $s=1,2,3$,
results in a linear recurrence for $Z_r(n,m)$ with~\eqref{Equ:ZRec} on the left-hand side and
$$
u(n,m)=\frac16\big(u_1(n,m)+u_2(n,m)+u_3(n,m)\big)
$$
on the right-hand side which holds for all $0\leq m\leq n-2$.
It remains to show that the inhomogeneous part evaluates to zero, 
$u(n,m)=0$ for $0\leq m\leq n-2$.
As indicated earlier, the expression $u(n,m)$ is composed by
\begin{itemize}
	\item the infinite sums~\eqref{Equ:InfiniteSums} with~\eqref{Equ:cDef};
	\item finite definite sums like those given in~\eqref{Equ:FiniteSums}.
\end{itemize}
A verification for all $n-2\geq m\geq0$ looks rather challenging. However, using the toolbox of~\texttt{Sigma}, this task can be accomplished automatically in 16~minutes of calculation time.
First, we treat the infinite sums by merging them to one big infinite sum and then compute a linear recurrence for it, which happens to be completely solvable in terms of indefinite nested sums.
This reduces all the infinite sums to indefinite nested sums.
The finite definite sums are a tougher nut to crack.
Internally, all sums (including~\eqref{Equ:FiniteSums}) are first considered as indefinite nested versions (with a common upper bound, say $a$). Then a finite subset of the sums arising is calculated with the command \texttt{SigmaReduce} such that there are no dependences among them and such that all the remaining sums can be represented in terms of these independent sums. 
It turns out that all sums (with $a$ now replaced by the `synchronized' upper bound $n-3$) cancel and only one definite sum remains. Activating the package \texttt{EvaluateMultiSums}~\cite{Schneider:13a} (that combines automatically the available summation tools of \texttt{Sigma}) this remaining sum simplifies to
\begin{align*}
&
\sum_{i=1}^{n-3}\frac{(-1)^i}i\binom{n}{i}\binom{2n-m+i}{n}
\\ &\quad
=\frac{(-1)^n \binom{3n-m-2}{n-2}}{2 (n-2) (n-1)^2 n^2}
\big(
-4 m
-4 m^2
+12 n+30 m n+6 m^2 n
-54 n^2
-43 m n^2-7 m^2 n^2
\\ &\quad\qquad
+70 n^3
+40 m n^3+4 m^2 n^3
-54 n^4
-19 m n^4
-m^2 n^4
+22 n^5
+4 m n^5
-4 n^6
\big)
\\ &\quad\;
-\binom{2 n-m}{n}\big(S_1(2n-m)+S_1(n)-S_1(n-m)\big).
\end{align*}
In a nutshell, $u(n,m)$ can be reduced to an expression given purely in terms of indefinite nested sums, which after further simplifications collapses to zero. This shows that not only the left-hand side but also the right-hand side of~\eqref{hyp1} satisfies the same recurrence~\eqref{Equ:ZRec}. The verification of this fact took in total 43 minutes.

\section{Dealing with the initial values}
\label{sec5}

In order to verify~\eqref{hyp1}, it remains to show~\eqref{Equ:Initial0} and~\eqref{Equ:Initial1}. For~\eqref{Equ:Initial0} we proceed as follows. First, we compute for $Z_l(n,0)$ the recurrence
\begin{multline}\label{Equ:RecIni0LHS}
-16 (2 n+1)^4 Z_l(n,0)
-(n+1)^4 Z_l(n+1,0)
\\
=-\frac{(-1)^n n!^8 (1+2 n)_{2 n}(1+4 n) \big(
	831+5265 n+12601 n^2+13499 n^3+5460 n^4\big)}{48 (2 n+1)!^5}.
\end{multline}
Internally, we follow the strategy in Section~\ref{Subsec:Z_l}: we use the representation from~\eqref{ZW1Ws} to get
$$
Z_l(n,0)=-\frac1{3}\big(W_1(n,0)+W_2(n,0)\big)
$$ 
and, for $W_1(n,0)$ and $W_2(n,0)$, compute two recurrences, where both have the \emph{same} homogeneous part.
Thus adding the inhomogeneous parts and simplifying the result further leads to~\eqref{Equ:RecIni0LHS}.
Solving this recurrence leads, for any $n\ge0$, to the closed form
\begin{align}
Z_l(n,0)
&=\frac{(-1)^n}{30720}
\big(
105 U_9(n)
+955 U_8(n)
+3095 U_7(n)
+2045 U_6(n)
\nonumber\\ &\;\quad
-12140 U_5(n)
-27300 U_4(n)
+12288 \zeta(2)^2
\big) \binom{2 n}{n}^4
\nonumber\\ &\;
+\frac{(-1)^n(4n+1)(5460 n^4+13499 n^3+12601 n^2+5265 n+831) \binom{4 n}{2 n}}{768(2n+1)^9\binom{2 n}{n}^4}
\label{Equ:Zln=0ClosedForm}
\end{align}
in terms of indefinite nested sums
\begin{equation}\label{eq:U}
U_k(n)=\sum_{i=0}^n\frac{\binom{4i}{2i}}{(2i+1)^k\binom{2i}{i}^8}
\quad\text{with}\; k=1,2,\dotsc.
\end{equation}
Similarly to Section~\ref{Subsec:Z_r}, we use the sum representation in~\eqref{Equ:ZrSumExpr} with $m=0$ encoded by $A_0(n,0)+\dots+A_3(n,0)$ to compute the recurrence
\begin{equation}\label{Equ:Zrm0}
\begin{split}
&
16 (n+1)^3 (2 n+1)^4 (4 n+3) (4 n+5) (5460 n^4+35339 n^3+85858 n^2+92804 n+37656) Z_r(n,0)
\\ &\;
+(357913920 n^{13}
+5716680688 n^{12}
+41762423804 n^{11}
+184637211081 n^{10}
\\ &\;\quad
+550778114541 n^9
+1169740743051 n^8
+1818232366245 n^7
+2092705983417 n^6
\\ &\;\quad
+1782121652067 n^5
+1108272850929 n^4
+488951050619 n^3
\\ &\;\quad
+144869028586 n^2
+25833166356 n
+2094206184) Z_r(n+1,0)
\\ &\;
+8 (n+2)^4 (2 n+3)^5  (5460 n^4+13499 n^3+12601 n^2+5265 n+831) Z_r(n+2,0)
=0
\end{split}
\end{equation}
which holds true for all $n\geq0$.
Furthermore, we verify that $Z_l(n,0)$ is also a solution of this recurrence by plugging its representation~\eqref{Equ:Zln=0ClosedForm} into the recurrence and checking that the expression simplifies to zero.
Finally, we verify that the first two initial values of $Z_l(n,0)$ and $Z_r(n,0)$ agree: 
\begin{equation*}
Z_l(0,0)=Z_r(0,0)=\frac25\zeta(2)^2,
\quad
Z_l(1,0)=Z_r(1,0)=\frac{277}{16} - \frac{32}{5}\zeta(2)^2;
\end{equation*}
to determine these evaluations again \texttt{Sigma} has been utilized. Together with the fact that the leading coefficient in~\eqref{Equ:Zrm0} is nonzero for all $n\geq0$, this implies that~\eqref{Equ:Initial0} holds. 

To verify~\eqref{Equ:Initial1}, we repeat the same game for $Z_l(n,1)$ and $Z_r(n,1)$: namely, we find the closed form representation
\begin{align}
Z_l(n,1)
&=
\frac{3 n(-1)^n}{40960}\big(
105 U_9(n)
+955 U_8(n)
+3095 U_7(n)
+2045 U_6(n)
\nonumber\\ &\;\quad
-12140 U_5(n)
-27300 U_4(n)
+12288 \zeta(2)^2
\big) \binom{2 n}{n}^4
\nonumber\\ &\;
-\frac{(-1)^n\binom{4 n}{2 n}}{1024 n^3 (2 n+1)^9\binom{2 n}{n}^4}\,
(16 n^9 + 116544 n^8 + 398115 n^7 + 587145 n^6
\nonumber\\ &\;\quad
+ 490329 n^5 + 255555 n^4 + 86016 n^3 + 18432 n^2 + 2304 n + 128)
\label{Equ:Zln=1ClosedForm}
\end{align}
valid for all $n\geq1$. In addition, we compute a recurrence of order 2 for $Z_r(n,1)$
and, as above, verify that $Z_l(n,1)$ is also its solution (by plugging in the representation~\eqref{Equ:Zln=1ClosedForm}).
Together with the initial values 
\begin{equation*}
Z_l(1,1)=Z_r(1,1)=-13 +\frac{24}{5}\zeta(2)^2,
\quad
Z_l(2,1)=Z_r(2,1)=\frac{4090247}{1944} - \frac{3888}{5}\zeta(2)^2
\end{equation*}
this implies that \eqref{Equ:Initial1} holds as well and completes the proof of~\eqref{hyp1}.
We note that the verification of each initial value problem, \eqref{Equ:Initial0} and~\eqref{Equ:Initial1}, took about 25 minutes.

\section{Summary}
\label{sec6}

Summarizing, the full proof of~\eqref{hyp1} took in total around 2 hours (excluding all the human trials and errors to find the tailored paths described above, and days to physically write this paper).

The initial values~\eqref{Equ:Zln=0ClosedForm} and~\eqref{Equ:Zln=1ClosedForm} are given through $2\zeta(2)^2/5=\zeta(4)$, hypergeometric products and the indefinite nested sums \eqref{eq:U} with $k=4,5,6,7,8,9$.
Thus, feeding the recurrence~\eqref{Equ:ZRec} with all this stuff we get the following corollary.

\begin{theorem}
For any $n\ge m\ge0$, both sides of $Z_l(n,m)=Z_r(n,m)$ can be expressed \textup(and computed in linear time\textup) in terms of $\zeta(4)$ and $U_4(n),\dots,U_9(n)$ in~\eqref{eq:U}.
\end{theorem}

The project \cite{MZ20} implicitly suggests that there can be further\,---\,more general(!)\,---\,forms of \eqref{hyp1}, with more than two independent parameters.
We have tried (unsuccessfully) to find some but cannot even figure out how to adopt \eqref{hyp1} to the case $m>n$.

\medskip
\noindent
\textbf{Acknowledgements.}
This project commenced during the joint visit of the authors in the Max Planck Institute for Mathematics (Bonn) in 2007
and went on during the second author's visit in the Research Institute for Symbolic Computation (Linz) in February 2020.
We thank the staff of these institutes for providing such excellent conditions for research.

\end{document}